\UseRawInputEncoding
\documentclass[12pt,reqno]{amsart}
\usepackage{amsmath,amsfonts,amsthm,amsopn,amssymb}
\usepackage{cite,marginnote}
\usepackage{color,enumitem,graphicx}
\usepackage[colorlinks=true,urlcolor=blue,
citecolor=red,linkcolor=blue,linktocpage,pdfpagelabels,
bookmarksnumbered,bookmarksopen]{hyperref}
\usepackage[english]{babel}

\usepackage[left=2.9cm,right=2.9cm,top=2.8cm,bottom=2.8cm]{geometry}
\usepackage[hyperpageref]{backref}

\numberwithin{equation}{section}

\makeindex
\makeindex
\newtheorem{theorem}{Theorem}[section]
\newtheorem{corollary}[theorem]{Corollary}

\newtheorem{lemma}[theorem]{Lemma}
\newtheorem{proposition}[theorem]{Proposition}

\theoremstyle{definition}
\newtheorem{definition}[theorem]{Definition}
\newtheorem{remark}[theorem]{Remark}

\newcommand{\RN}{\mathbb R^N}

\newcommand{\s}{\section}

\newcommand{\R}{\mathbb R}

\newcommand{\lab}{\label}
\newcommand{\bt}{\begin{theorem}}
\newcommand{\et}{\end{theorem}}
\newcommand{\bl}{\begin{lemma}}
\newcommand{\el}{\end{lemma}}
\newcommand{\bd}{\begin{definition}}
\newcommand{\ed}{\end{definition}}
\newcommand{\bc}{\begin{corollary}}
\newcommand{\ec}{\end{corollary}}
\newcommand{\bpr}{\begin{proof}}
\newcommand{\epr}{\end{proof}}
\newcommand{\bx}{\begin{example}}
\newcommand{\ex}{\end{example}}
\newcommand{\bi}{\begin{exercise}}
\newcommand{\ei}{\end{exercise}}
\newcommand{\bo}{\begin{proposition}}
\newcommand{\eo}{\end{proposition}}
\newcommand{\br}{\begin{remark}}
\newcommand{\er}{\end{remark}}
\newcommand{\beq}{\begin{equation}}
\newcommand{\eeq}{\end{equation}}
\newcommand{\ba}{\begin{align}}
\newcommand{\ea}{\end{align}}
\newcommand{\bn}{\begin{enumerate}}
\newcommand{\en}{\end{enumerate}}
\newcommand{\bg}{\begin{align*}}
\newcommand{\bcs}{\begin{cases}}
\newcommand{\ecs}{\end{cases}}

\newcommand{\bean}{\begin{eqnarray*}}
\newcommand{\eean}{\end{eqnarray*}}

\def\N{\mathbb{N}}

\def\R{\mathbb{R}}

\def\bd{\mathrm{bd}\,}

\title[Normalized solution to at least mass critical problem]{On the Kirchhoff equation with prescribed mass and general nonlinearities}

\author[X.~Y.~Zeng]{Xiaoyu Zeng}
\author[J.~J.~Zhang]{Jianjun Zhang}
\author[Y.~M.~Zhang]{Yimin Zhang}
\author[X.~X.~Zhong]{Xuexiu Zhong}

\address[X.~Y.~Zeng]{\newline\indent ~Center for Mathematical Sciences
\newline\indent
Wuhan University of Technology,
\newline\indent
Wuhan, 430070, China.}
\email{\href{mailto:xyzeng@whut.edu.cn}{xyzeng@whut.edu.cn}}

\address[J.~J.~Zhang]{\newline\indent College of Mathematics and Statistics
\newline\indent
Chongqing Jiaotong University
\newline\indent
Chongqing 400074, PR China}
\email{\href{mailto:zhangjianjun09@tsinghua.org.cn}{zhangjianjun09@tsinghua.org.cn}}

\address[Y.~M.~Zhang]{\newline\indent Center for Mathematical Sciences
\newline\indent
Wuhan University of Technology
\newline\indent
Wuhan, 430070,PR China}
\email{\href{mailto:zhangym802@126.com}{zhangym802@126.com}}

\address[X.~X.~Zhong]{\newline\indent South China Research Center for Applied Mathematics and Interdisciplinary Studies
\newline\indent
South China Normal University
\newline\indent
Guangzhou 510631, PR China}
\email{\href{mailto:zhongxuexiu1989@163.com}{zhongxuexiu1989@163.com}}

\thanks{$^*$Corresponding author: Xuexiu Zhong}

\thanks{The first author is supported by NSFC (No.12171379). The second author is supported by the NSFC (No.11871123). The third author is partially supported by the NSFC (No.11771127) and the Fundamental Research Funds for the Central Universities (WUT: 2020IB011). The corresponding author is partially supported by the NSFC (No.12271184), Guangdong Basic and Applied Basic Research Foundation (2021A1515010034),Guangzhou Basic and Applied Basic Research Foundation(202102020225).}

\subjclass[2000]{}

\begin{document}

\begin{abstract}
In the present paper, we apply a global branch approach to study the existence, non-existence and multiplicity of positive normalized solutions $(\lambda_c, u_c)\in \mathbb{R}\times H^1(\mathbb{R}^N)$ to the following Kirchhoff problem
$$
-\left(a+b\int_{\mathbb{R}^N}|\nabla u|^2dx\right)\Delta u+\lambda u=g(u)~\hbox{in}~\mathbb{R}^N,\;N\geq 1
$$
satisfying the normalization constraint
$
\displaystyle\int_{\mathbb{R}^N}u^2=c,
$
which appears in free vibrations of elastic strings.  The parameters $a,b>0$ are prescribed as well as the mass $c>0$. Due to the presence of the non-local term $b\int_{\mathbb{R}^N}|\nabla u|^2dx \Delta u$, such problems lack the mountain pass geometry in the higher dimension case $N\geq 5$. Our result seems to be the first attempt in this aspect.
\vskip0.1in
{\large \it Dedicated to Prof. Vicen\c{t}iu D. R\u{a}dulescu on the occasion of his 65th birthday.}
\vskip0.1in
\end{abstract}

\maketitle

\s{Introduction}
\renewcommand{\theequation}{1.\arabic{equation}}

In the present paper, we are concerned about the following normalization problem: to find one couple $(\lambda, u)\in \R\times H^1(\R^N)$ solving
\beq\lab{1.1}
-\left(a+b\int_{\R^N}|\nabla u|^2dx\right)\Delta u+\lambda u=g(u)~\hbox{in}~\R^N,\;N\geq 1
\eeq
with
\beq\lab{eq:20210816-constraint}
\int_{\R^N}|u|^2 dx=c,
\eeq
where $a,b,c>0$ are prescribed.

The Kirchhoff equation \eqref{1.1} is an extension of the classical D'Alembert's wave equations for free vibration of elastic strings. In order to describe the transversal free vibrations of a clamped string in which the dependence of the tension on the deformation cannot be neglected,  Kirchhoff \cite{Kirchhoff} firstly introduced such an equation in dimension 1, without forcing term and with Dirichlet boundary conditions. From the mathematical point of view, problem \eqref{1.1} is nonlocal since the appearance
of the term $\int_{\R^N}|\nabla u|^2dx \Delta u$. It is worth mentioning that a functional analysis approach was introduced by J.~L. Lions in the pioneer work \cite{Lions1978}. After then, \eqref{1.1} has attracted extensive attentions.

For the {\it fixed frequency problem}, i.e., $\lambda$ is prior given, the solutions to \eqref{1.1} correspond to the critical points of the associated energy functional $\Phi_\lambda: H^1(\R^N)\mapsto \R$, defined by
$$
\Phi_\lambda[u]:=\frac{a}{2}\int_{\mathbb{R}^N}|\nabla u|^2dx+\frac{b}{4}\left(\int_{\mathbb{R}^N}|\nabla u|^2dx\right)^2+\frac{\lambda}{2}\int_{\R^N}|u|^2 dx-\int_{\mathbb{R}^N}G(u)dx,
$$
where $G(s):=\int_0^s g(t)dt$.

From the viewpoint of calculus of variations, since the order of the Kirchhoff type
nonlocal term  $-b\left(\int_{\R^N}|\nabla u|^2dx\right)\Delta u$ in the corresponding functional is $4$, the following assumptions are classical.
 \begin{itemize}
\item[$\mathbf{\mathbf{(g_1)}}$]$\lim_{|s|\rightarrow +\infty}\frac{G(s)}{|s|^4}=\infty, ~\hbox{where}~G(s):=\int_0^s g(t)dt$.

\item[$\mathbf{\mathbf{(g_2)}}$] $\frac{g(s)}{|s|^3}~\hbox{increases strictly for}~t\in \R\backslash\{0\}$.

\item[$\mathbf{\mathbf{(g_3)}}$]There exists some $\eta>4$ such that $0<\eta G(s)\leq g(s)s$ for $s\neq 0$.
\end{itemize}

 One usually assume $\mathbf{\mathbf{(g_1)}}$ to apply the mountain pass theorem or  $\mathbf{\mathbf{(g_2)}}$ to apply the Nehari manifold method.
Furthermore, to guarantee the boundedness of a Palais-Smale sequence for $\Phi_\lambda$, Ambrosetti-Rabinowitz type condition $\mathbf{\mathbf{(g_3)}}$ is usually adopted.

For the case of $N\geq 3$, the critical Sobolev embedding  states out that the nonlinearity cannot exceed the polynomial of degree $2^*=\frac{2N}{N-2}$. We remark that the Sobolev exponent $2^*=4$ in the case of $N=4$ and $2^*<4$ for $N\geq 5$. In particular, for the case of $N\geq 5$, if the nonlinearity has a growth at most Sobolev critical, which mean that the nonlinearity has a growth order less than 4, then the action functional $\Phi_\lambda$ with $\lambda$ is bounded from below and coercive in $H^1(\R^N)$. And thus one can apply a minimization argument to find a critical point. However, this is not the application scope of the mountain pass theorem. To get a mountain pass solution, almost all literatures are only concerned with \eqref{1.1} for $N\leq 3$, see \cite{Alves2005,Azzollini2015,Azzollini2012,Zhang2020,Deng2015,Guo2015,He2012,He2014,He2016,Wang2012,Wang2020,Shuai2015,Naimen2014,Figueiredo2014,Lei2015,Li2014,Li2012} and the references therein.
The existence of positive solutions are studied
in \cite{Alves2005,Guo2015,Figueiredo2014,Li2014,Li2012,Naimen2014}. And the existence of ground state solutions, see \cite{Azzollini2015,Azzollini2012,Guo2015,Li2014}. The multiplicity results of problems related \eqref{1.1}, one can found in \cite{Zhang2020,Wang2012,Lei2015,He2012,Chen2013}.
The asymptotic behavior of the solutions is also a research interest of scholars.
In \cite{He2012,Figueiredo2014,Wang2012,Wang2020}, the authors considered \eqref{1.1} with very general nonlinearities, even involving the Sobolev critical exponent, and studied the existence, multiplicity and concentration behavior of positive solutions.  For the bounded states, we refer to \cite{He2014,He2016}. And the sign-changing solutions can be seen in \cite{Deng2015,Shuai2015}.

In \cite{Azzollini2015}, Azzollini  gave a characterization of the solution to the Kirchhoff equation \eqref{1.1} as follows.

\vspace*{4pt}\noindent
{\bf Theorem A.}(cf.\cite[Theorem 1.1]{Azzollini2015}) \emph{Fix $N\geq 1, a>0,b>0,\lambda>0$. $u\in C^2(\R^N)\cap H^{1}(\R^N)$ is a solution to \eqref{1.1}
 if and only if there exists $v\in C^2(\R^N)\cap H^{1}(\R^N)$ solving}
 \beq\lab{eq:20211029-wbe1}
 -\Delta v+\lambda v=g(v)~\hbox{in}~\R^N,
 \eeq
\emph{and $t>0$ such that}
\beq\lab{eq:20211119-e4}
t^2 a+t^{4-N}b\|\nabla v\|_2^2=1~\hbox{and}~u(x)=v(tx).
\eeq
\qed

Azzollini's correspondence  provides a way to search solutions of \eqref{1.1} for higher dimensional cases. Indeed, for $N\geq 4$ and any solution $v\in H^1(\R^N)$ to \eqref{eq:20211029-wbe1}, the algebraic equation  in \eqref{eq:20211119-e4} always possesses a root $t\in \R^+$ provided $b$ small, and thus $u(x):=v(tx)$ is a solution to \eqref{1.1}.

We remark that the euqation
 $$
 -\Delta v=f(v)~\hbox{in}~\R^N, 0\not\equiv v\in H^1(\R^N)
$$
 has been well studied in \cite{Berestycki1983,Berestycki1983a,Berestycki1983b,Jeanjean2003,Jeanjean2003a} under the following very general assumptions
\begin{itemize}
 \item[$\mathbf{\mathbf{(f_1)}}$] $f\in C(\R,\R)$ and is of odd.
 \item[$\mathbf{\mathbf{(f_2)}}$]$-\infty<\liminf_{s\rightarrow 0}\frac{f(s)}{s}\leq \limsup_{t\rightarrow 0}\frac{f(s)}{s}=-\omega<0$ for $N\geq 3$.\\
     $\lim_{t\rightarrow 0}\frac{f(s)}{s}=-\omega\in (-\infty,0)$ for $N=1,2$.
 \item[$\mathbf{\mathbf{(f_3)}}$]When $N\geq 3$, $\lim_{s\rightarrow \infty} \frac{f(s)}{|s|^{\frac{N+2}{N-2}}}=0$.\\
     When $N=2$, for any $\alpha>0$, $\lim_{s\rightarrow \infty} \frac{f(s)}{e^{\alpha s^2}}=0$.
 \item[$\mathbf{\mathbf{(f_4)}}$]Let $F(s):=\int_0^s f(t)dt$. When $N\geq 2$, there exists $\xi>0$ such that $F(\xi)>0$. When $N=1$, there exist $\xi>0$ such that
     $$F(s)<0, \forall s\in (0,\xi), \;F(\xi)=0~\hbox{and}~f(\xi)>0.$$
 \end{itemize}
Then basing on {\bf Theorem A}, also under the assumptions $\mathbf{\mathbf{(f_1)}}$--$\mathbf{\mathbf{(f_4)}}$, Lu \cite{Lu2017} generalized the results to the Kirchhoff equation
$$
-\left(a+b\int_{\R^N}|\nabla u|^2dx\right)\Delta u=f(u)~\hbox{in}~\R^N,\;N\geq 1.
$$
Even for the high-dimensional case with nonlinearity involving Sobolev critical exponent, the Azzollini's correspondence is powerful, see \cite{Xie2022}. We emphasize that both \cite{Lu2017} and \cite{Xie2022}, the parameter $b$ is assumed small proper for the high dimension case, aim to guarantee a positive root for the algebra equation \eqref{eq:20211119-e4}.

There is also another way to study problem \eqref{1.1}, by searching for solutions satisfying the mass constraint \eqref{eq:20210816-constraint}.
In such a point of view, the mass $c>0$ is prescribed, while the frequency $\lambda$ is unknown and will come out as a Lagrange multiplier. Hence, we call it {\it fixed mass problem} and the solution $u$ is called normalized solution. Nowadays, the study of normalized problem has attracted a lot of mathematical researchers' concerns, due to the physic reasons. Naturally, normalized solutions of \eqref{1.1} can be searched as critical points of $I[u]$ constrained on $S_c$, where
$$
I[u]:=\Phi_0[u]=\frac{a}{2}\int_{\mathbb{R}^N}|\nabla u|^2dx+\frac{b}{4}\left(\int_{\mathbb{R}^N}|\nabla u|^2dx\right)^2-\int_{\mathbb{R}^N}G(u)dx,
$$
and
$$
S_c:=\left\{u\in H^1(\R^N): \|u\|_2^2=c\right\}.
$$

 As far as we know, for the study of normalized solution, the first work for equation \eqref{1.1} with $g(u)=|u|^{p-2}u$ is due to Ye. And it has been studied very well in a sequence of his work \cite{Ye2015,Ye2015a,Ye2016}. The following minimization problem is considered:
\begin{equation}\label{eq1.2}
E_c:=\inf_{u\in S_c}I[u].
\end{equation}

Ye proved the existence and non-existence of minimizer by a scaling technique and applying the concentration-compactness principle, provided $p<2+\frac{8}{N}$. When $2+\frac{8}{N}< p<2^*$, he can get a  mountain pass type solution.  If $p=2+\frac{8}{N}$, the asymptotic behavior of critical points of $I[u]$ on $S_c$ is also studied as $c\uparrow \bar{c}$.
By scaling technique and energy estimate, Zeng and Zhang \cite{Zeng2017} improved the results of \cite{Ye2015}. Precisely, they obtained the existence and uniqueness of minimizers of (\ref{eq1.2}) with $2<p<2+\frac{8}{N}$, and the existence and uniqueness of the mountain pass type critical points on the $L^2$ normalized manifold for $2+\frac{8}{N}<p<2^*$ or $p=2+\frac{8}{N}$ and $c>\bar{c}$ for some critical mass $\bar{c}>0$.

Using minimax procedure, Luo and Wang \cite{Luo2017} obtained the multiplicity of critical points of $I(u)$ on the $L^2$ normalized manifold for $2+\frac{8}{N}<p<2^*$ and any $c>0$. Recently, He, Lv, the third author and the fourth author \cite{HeLvZhangZhong2021} consider the general nonlinearity $g$ of mass super critical case, they study the existence of ground state normalized solutions for any given $c>0$ via Pohozaev manifold constraint method and also study the asymptotic behavior of these solutions as $c\rightarrow 0^+$ as well as $c\rightarrow+\infty$, via a sequence of blow up arguments. If $g$ has Sobolev critical exponent term, the existence of critical points of $I(u)$ constrained on the $L^2$-sphere was obtained in \cite{Li2021,Zhang2021}. Other results about existence of normalized solutions for equation \eqref{1.1} with potential function, can be seen in \cite{Guo2018,Li2019,Meng2022,Zhu2021}. We emphasize that results in \cite{Ye2015,Ye2015a,Ye2016,Zeng2017,Luo2017,HeLvZhangZhong2021,Li2021,Zhang2021,Guo2018,Li2019,Meng2022,Zhu2021} are only considered in dimension $N\leq 3$, and can not be generalized to high dimension case.

We note that when $b=0$, the normalized solution problem
 \beq\lab{eq:20211119-e7}
 \begin{cases}
 -\Delta v+\lambda v=g(v)~\hbox{in}~\R^N,\\
 \int_{\R^N}|v|^2 dx=c
 \end{cases}
 \eeq
 also has been studied in a large number of literatures. Let us just mention here, see\cite{Stuart1981,Stuart1982,Lions1984a,Lions1984b,Shibata2014, JeanLu2019, Jeanjean-Lu-2022, Stefanov19,Jeanjean1997,BartschValeriola2013, BartschSoave2017, BieganowskiMederski, JeanLu2020a, ikomatanaka2019, Soave-1,JeanjeanLe2022, Soave2020, Wei-Wu2022} and the references therein.

\br\lab{remark:20211119-r1}
We remark that unlike \cite{Lu2017,Xie2022},  the Azzollini's correspondence is not applied in the study of normalized solutions of \eqref{eq:20211119-e7} via a variational method, even for small $b$. Since the mass constraint condition  \eqref{eq:20210816-constraint} is not kept in the transformation when doing a Azzollini's correspondence.
 \er

Recently, Jeanjean, the second author and the fourth author \cite{JeanZhangZhong2021} develop a novel approach to the normalized solutions of Schrodinger equation \eqref{eq:20211119-e7}.
It is worth mentioning that the classical approach to the normalized solution problem is variational on some suitable manifolds, which heavily depends on the geometry of the energy functional. Unlike the classical approach by variational, such a global branch approach does not depend on the geometry of the energy functional. So that, they can handle the nonlinearities $g$ in a unified way, which are either mass subcritical, mass critical or mass supercritical.

For the convenience of discussion, it is time to put out the conditions involved in  this work. For the nonlinear term $g(s)$, we  always assume the  conditions as follows.
\begin{itemize}
\item[$\mathbf{\mathbf{(G_1)}}$]$g \in C^1([0, + \infty))$, $g(s)>0$ for $s> 0$.
\item[$\mathbf{\mathbf{(G_2)}}$] There exists some $(\alpha,\beta)\in \R_+^2$ satisfying
$$\begin{cases}
2<\alpha,\beta<2^*:=\frac{2N}{N-2}\quad&\hbox{if}~N\geq 3,\\
2<\alpha, \beta<2^*:=+\infty \quad&\hbox{if}~N=1,2,
\end{cases}$$
such that
$$\lim_{s\rightarrow 0^+}\frac{g'(s)}{s^{\alpha-2}}=\mu_1(\alpha-1)>0~\hbox{and}~\lim_{s\rightarrow +\infty}\frac{g'(s)}{s^{\beta-2}}=\mu_2(\beta-1)>0.$$
\item[$\mathbf{\mathbf{(G_3)}}$] $-\Delta u=g(u)$ has no positive radial decreasing classical solution in $\mathbb{R}^N$.
\end{itemize}

We also denote $$I_1=\left(2, 2+\frac{4}{N}\right), I_2=\left\{2+\frac{4}{N}\right\}, I_3=\left(2+\frac{4}{N}, 2^*\right)$$
and
$$J_1=\left(2, 2+\frac{8}{N}\right), J_2=\left\{2+\frac{8}{N}\right\}, J_3=\left(2+\frac{8}{N}, 2^*\right).$$
For simplicity, we write
$$
K_{ij}:=I_i\times J_j, \quad i,j\in \{1,2,3\}.
$$

Let $U \in C_{r,0}(\RN)$ be the unique positive solution to
\beq\lab{eq:20210820-e6L}
-\Delta U+U=\mu_1 U^{\alpha-1}~\hbox{in}~\RN,
\eeq
and $V \in C_{r,0}(\RN)$ be the unique positive solution to
\beq\lab{eq:20210820-e6l}
-\Delta V+V=\mu_2 V^{\beta-1}~\hbox{in}~\RN,
\eeq
where $C_{r,0}(\RN)$ denotes the space of continuous radial functions vanishing at $\infty$.
Here comes  our first main result for $N\leq 3$.
\bt\lab{th:20211109-th1}({\bf The case of $N\leq3$})
Let $a>0,b>0, N\leq 3$ and $\mathbf{\mathbf{(G_1)}}$--$\mathbf{\mathbf{(G_3)}}$ hold. Then we have the following conclusions.
\begin{itemize}
\item[(1)]If $(\alpha,\beta)\in K_{11}$,
    then for any give $c>0$, \eqref{1.1}-\eqref{eq:20210816-constraint} possesses a positive normalized solution $(\lambda, u_\lambda)\in (0, + \infty) \times H_{rad}^{1}(\R^N)$.
\item[(2)]If $(\alpha,\beta)\in K_{21}$, then \eqref{1.1}-\eqref{eq:20210816-constraint} possesses a positive normalized solution $(\lambda, u_\lambda)\in (0, + \infty) \times H_{rad}^{1}(\R^N)$ provided $c\in (a^N\|U\|_2^2,+\infty)$. Furthermore, there exists some $0<c^*\leq a^N\|U\|_2^2$ such that \eqref{1.1}-\eqref{eq:20210816-constraint} have no positive normalized solution if $c\in (0,c^*)$.
\item[(3)]If $(\alpha,\beta)\in K_{31}$, there exists some $0<m\leq c^*<+\infty$ such that for any $c\in (c^*,+\infty)$, \eqref{1.1}-\eqref{eq:20210816-constraint}  have at least two distinct positive normalized solutions $(\lambda_i, u_{\lambda_i}) \in (0, + \infty) \times H_{rad}^{1}(\R^N),i=1,2$, while \eqref{1.1}-\eqref{eq:20210816-constraint} have no positive normalized solution provided $c<m$.
\item[(4)]If $(\alpha,\beta)\in K_{12}$, then for any $c\in \left(0,b^{\frac{N}{4-N}}\left(\|\nabla V\|_2^2\right)^{\frac{N}{4-N}} \|V\|_2^2\right)$, \eqref{1.1}-\eqref{eq:20210816-constraint}  have at least one positive normalized solutions $(\lambda, u_{\lambda}) \in (0, + \infty) \times H_{rad}^{1}(\R^N) $. Furthermore, there exists some $c^*>b^{\frac{N}{4-N}}\left(\|\nabla V\|_2^2\right)^{\frac{N}{4-N}} \|V\|_2^2$ such that \eqref{1.1}-\eqref{eq:20210816-constraint} have no positive normalized solution provided $c>c^*$.
\item[(5)]If $(\alpha,\beta)\in K_{32}$, there exists some $0<m<b^{\frac{N}{4-N}}\left(\|\nabla V\|_2^2\right)^{\frac{N}{4-N}} \|V\|_2^2$ such that for any $c\in \left(b^{\frac{N}{4-N}}\left(\|\nabla V\|_2^2\right)^{\frac{N}{4-N}} \|V\|_2^2, +\infty\right)$, \eqref{1.1}-\eqref{eq:20210816-constraint}  have at least one positive normalized solution $(\lambda, u_{\lambda}) \in (0, + \infty) \times H_{rad}^{1}(\R^N) $, while \eqref{1.1}-\eqref{eq:20210816-constraint} have no positive normalized solution provided $c<m$.
\item[(6)]If $(\alpha,\beta)\in K_{13}$, there exists some $0<c^*\leq M<+\infty$ such that for any $c\in (0,c^*)$, \eqref{1.1}-\eqref{eq:20210816-constraint}  have at least two distinct positive normalized solutions $(\lambda_i, u_{\lambda_i}) \in (0, + \infty) \times H_{rad}^{1}(\R^N),i=1,2 $, while \eqref{1.1}-\eqref{eq:20210816-constraint} have no positive normalized solution provided $c>M$.
\item[(7)]If $(\alpha,\beta)\in K_{23}$, then for any $c\in (0,a^N\|U\|_2^2)$, \eqref{1.1}-\eqref{eq:20210816-constraint}  have at least one positive normalized solutions $(\lambda, u_{\lambda}) \in (0, + \infty) \times H_{rad}^{1}(\R^N) $. Furthermore, there exists some $ M\geq a^N\|U\|_2^2$ such that, \eqref{1.1}-\eqref{eq:20210816-constraint} have no positive normalized solution provided $c>M$.
\item[(8)]If $(\alpha,\beta)\in K_{33}$,
    then for any given $c>0$, \eqref{1.1}-\eqref{eq:20210816-constraint} possesses a positive normalized solution $(\lambda, u_\lambda)\in (0, + \infty) \times H_{rad}^{1}(\R^N)$.
\item[(9)]If $(\alpha,\beta)\in K_{22}$, we put $$c_1=\min\left\{a^N\|U\|_2^2, b^{\frac{N}{4-N}}\left(\|\nabla V\|_2^2\right)^{\frac{N}{4-N}} \|V\|_2^2\right\}$$
    and
    $$c_2=\max\left\{a^N\|U\|_2^2, b^{\frac{N}{4-N}}\left(\|\nabla V\|_2^2\right)^{\frac{N}{4-N}} \|V\|_2^2\right\}.$$
    Then there exists some $m<c_1$ and $M>c_2$ such that \eqref{1.1}-\eqref{eq:20210816-constraint} have no positive normalized solution provided $c\in (0,m)\cup (M,+\infty)$. And if $(c_1,c_2)\neq \emptyset$, then \eqref{1.1}-\eqref{eq:20210816-constraint} possesses at one positive normalized solution $(\lambda, u_\lambda)\in (0, + \infty) \times H_{rad}^{1}(\R^N)$ provided $c\in (c_1,c_2)$.
\end{itemize}
\et

 For the high dimension case, we can also establish similar results about existence, non-existence and multiplicity. However, different from the case of $N\leq 3$, it is surprising that the results do not rely on $\beta$ for the case of $N\geq 4$. Precisely, we obtain the following results.

\bt\lab{thm:20211115-th1}({\bf The case of $N=4$})
Let $a>0,b>0, N=4$ and $\mathbf{\mathbf{(G_1)}}$--$\mathbf{\mathbf{(G_3)}}$ hold. Then we have the following conclusions.
\begin{itemize}
\item[(i)]If $2<\alpha <2+\frac{4}{N}$, then for any $c\in (0,+\infty)$, \eqref{1.1}-\eqref{eq:20210816-constraint}  have at least one positive normalized solutions $(\lambda, u_{\lambda}) \in (0, + \infty) \times H_{rad}^{1}(\R^N)$;
\item[(ii)] If $\alpha=2+\frac{4}{N}$, then for any $c\in (a^{\frac{N}{2}}\|U\|_2^2, +\infty)$, \eqref{1.1}-\eqref{eq:20210816-constraint}  have at least one positive normalized solution $(\lambda, u_{\lambda}) \in (0, + \infty) \times H_{rad}^{1}(\R^N)$;
\item[(iii)]If $2+\frac{4}{N}<\alpha<2^*$, then there exists some $c^*>0$ such that for any $c>c^*$, \eqref{1.1}-\eqref{eq:20210816-constraint}  have at least two distinct positive normalized solutions $(\lambda_i, u_{\lambda_i}) \in (0, + \infty) \times H_{rad}^{1}(\R^N),i=1,2$.
\end{itemize}
\et

\bt\lab{thm:20211115-th2}({\bf The case of $N\geq 5$})
Let $a>0,b>0, N\geq 5$ and $\mathbf{\mathbf{(G_1)}}$--$\mathbf{\mathbf{(G_3)}}$ hold. Then we have the following conclusions.
\begin{itemize}
\item[(i)]If $2<\alpha<2+\frac{4}{N}$, then for any $c\in (0,+\infty)$, \eqref{1.1}-\eqref{eq:20210816-constraint}  have at least one positive normalized solutions $(\lambda, u_{\lambda}) \in (0, + \infty) \times H_{rad}^{1}(\R^N)$;
\item[(ii)] If $\alpha=2+\frac{4}{N}$, then for any $c\in (a^{\frac{N}{2}}\|U\|_2^2, +\infty)$, \eqref{1.1}-\eqref{eq:20210816-constraint}  have at least one positive normalized solution $(\lambda, u_{\lambda}) \in (0, + \infty) \times H_{rad}^{1}(\R^N)$;
\item[(iii)]If $2+\frac{4}{N}<\alpha<2^*$, then there exists some
$$ 0<c^*\leq  (\frac{b(N-4)}{2a})^{-\frac{N}{N-2}}\|\nabla v(s^*, \cdot)\|_{2}^{-\frac{2N}{N-2}} \|v(s^*, \cdot)\|_2^2$$
such that for any $c>c^*$, \eqref{1.1}-\eqref{eq:20210816-constraint}  have at least two distinct positive normalized solutions $(\lambda_i, u_{\lambda_i}) \in (0, + \infty) \times H_{rad}^{1}(\R^N), i=1,2$. Here $s^*$ is given by Lemma \ref{lemma:20211115-l2} and $v(s,x)\in H^1(\R^N)$ is given by \eqref{eq:20230217-e1} below.
\end{itemize}
\et

For $N\geq 4$, we also obtain the following non-existence result for the positive normalized solution.
\bt\lab{thm:20211115-th3}
Let $a>0,b>0, N\geq 4$ and $\mathbf{\mathbf{(G_1)}}$--$\mathbf{\mathbf{(G_3)}}$ hold. If $\alpha\in I_2\cup I_3$, then there exists some $c_0>0$ such that \eqref{1.1}-\eqref{eq:20210816-constraint} has no positive solution provided $c<c_0$.
\et

When study the positive normalized solutions via a global branch approach,
the main ingredients is the study of the asymptotic behaviors of the positive solutions as $\lambda\rightarrow 0^+$ or $\lambda\rightarrow +\infty$ and the existence of an unbounded continuum of solutions in $(0, + \infty) \times H^1(\mathbb{R}^N)$. We remark that under the assumptions $\mathbf{\mathbf{(G_1)}}$--$\mathbf{\mathbf{(G_3)}}$, the authors in \cite{JeanZhangZhong2021} indeed manage to obtain a connected positive solutions set
\beq\lab{eq:20230217-e1}
\mathcal{\tilde{S}}:=\{(\lambda(s), v(s,x))\in \R^+\times H_{rad}^{1}(\R^N):s\in \R, \lim_{s\rightarrow -\infty}\lambda(s)=0,\lim_{s\rightarrow +\infty}\lambda(s)=+\infty\}
 \eeq
 to \eqref{eq:20211029-wbe1}. And they also make clear the behavior of positive solutions of \eqref{eq:20211029-wbe1} as $\lambda\rightarrow 0^+$ as well as $\lambda\rightarrow +\infty$, see \cite[Theorem 1.3]{JeanZhangZhong2021}.

In this work, our arguments are basing on a global branch approach and the Azzollini's correspondence, since preliminaries are different for the cases of $N\leq 3$, $N=4$ and $N\geq 5$. We prefer to deal with  them  individually. We shall establish some preliminaries for $N\leq 3$ and give the proof of  Theorem \ref{th:20211109-th1} in Section \ref{sec:20230217-s1}. Similarly, we focus on $N=4$ and give the proof of  Theorem \ref{thm:20211115-th1} in Section \ref{sec:20230217-s2}. In Section \ref{sec:20230217-s3}, we study the case of $N\geq 5$ and prove Theorem \ref{thm:20211115-th2}. Theorem \ref{thm:20211115-th3} is focus on the non-existence result, which will be proved in Section \ref{sec:20230217-s4}.

\section{The case of $N\leq 3$}\lab{sec:20230217-s1}
\renewcommand{\theequation}{2.\arabic{equation}}

\bl\lab{lemma:20211109-l1}
Let $a>0,b>0$ and $N\leq 3$. For any $v\in D_{0}^{1,2}(\R^N)$, there exists a unique $t_v>0$ such that $t_v$ solves $at^2+b\|\nabla v\|_2^2t^{4-N}-1=0$.
Furthermore, $t_v$ depends continuously on $v\in D_{0}^{1,2}(\R^N)$ and
$t_v=\frac{1}{\sqrt{a}}\big(1+o(1)\big)$
as $\|\nabla v\|_2\rightarrow 0$, while
$b\|\nabla v\|_2^2t_{v}^{4-N}=1+o(1)$ as $\|\nabla v\|_2\rightarrow +\infty$.
\el
\bpr
For any $v\in D_{0}^{1,2}(\R^N)$, we put $f(t):=at^2+b\|\nabla v\|_2^2t^{4-N}-1\in C^1([0,+\infty))$. Since $a>0,b>0,N\leq 3$, we have $f'(t)=2at+(4-N)b\|\nabla v\|_2^2t^{3-N}>0, t>0$, which implies that $f(t)$ increases strictly in $[0,+\infty)$. On the other hand, by $N\leq 3$, we see that $f(0)=-1<0$ and $\displaystyle \lim_{t\rightarrow +\infty}f(t)=+\infty$. Hence, there exists a unique $t_v>0$ such that $f(t_v)=0, f(t)<0$  for $t\in [0,t_v)$ and $f(t)>0$ for $t>t_v$.

Let $v_n\rightarrow v_0$ in  $D_{0}^{1,2}(\R^N)$, we have that $\|\nabla v_n\|_2^2\rightarrow \|\nabla v_0\|_2^2$ as $n\rightarrow \infty$. By the definition of $t_{v_n}$, it is easy to see that $\{t_{v_n}\}$ is bounded in $\R$. Up to a subsequence, still denoted by $t_{v_n}$, we may assume that $t_{v_n}\rightarrow \bar{t}\geq 0$. So
$$a \bar{t}^2+b \|\nabla v_0\|_2^2 \bar{t}^{4-N}-1=\lim_{n\rightarrow +\infty}\left[at_{v_n}^{2}+b\|\nabla v_n\|_2^2t_{v_n}^{4-N}-1\right]=0,$$
which implies that $t_{v_0}=\bar{t}$. That is, $t_{v_n}\rightarrow t_{v_0}$.

In particular, if $v_n\rightarrow 0$ in $D_{0}^{1,2}(\R^N)$, since $\{t_{v_n}\}$ is bounded, we have that $b\|\nabla v_n\|_2^2t_{v_n}^{4-N}=o(1)$ and thus $t_{v_n}=\frac{1}{\sqrt{a}}\big(1+o(1)\big)$.
If $\|\nabla v_n\|_2\rightarrow +\infty$, then by $b>0$ we have that $t_{v_n}\rightarrow 0$. Thus, $b \|\nabla v_n\|_2^2 t_{v_n}^{4-N}\rightarrow 1$ as $n\rightarrow +\infty$.
\epr

\bl\lab{lemma:20211109-l3}
Assume that $a>0,b>0$ and $N\leq 3$. Let $\tilde{\mathcal{S}}$ be the connected positive solutions  branch  given by \eqref{eq:20230217-e1}. For any $s>0$, let $t_{v(s,x)}$ be the unique number determined by Lemma \ref{lemma:20211109-l1} and denoted by $t_s$ for simplicity. Put $u(s,x):=v(s,t_s x)$. Then
\begin{itemize}
\item[(i)] $\mathcal{C}:=\left\{(\lambda(s), u(s,x)): s\in \R\right\}$ is a connected positive solution branch of \eqref{1.1}.
\item[(ii)]\beq\lab{eq:20211109-be1}
\lim_{s \rightarrow -\infty }\|\nabla u(s,\cdot)\|_2=0
\eeq
and
\beq\lab{eq:20211109-be2}
\lim_{s \rightarrow -\infty }\|u(s,\cdot)\|_2^2=\begin{cases}
    0\quad &\alpha<2+\frac{4}{N},\\
    a^{\frac{N}{2}}\|U\|_2^2& \alpha=2+\frac{4}{N},\\
    +\infty& \alpha>2+\frac{4}{N}.
    \end{cases}
\eeq
\item[(iii)]$\displaystyle\lim_{s \rightarrow +\infty }\|\nabla u(s,\cdot)\|_2=+\infty$
and
$$
\lim_{s \rightarrow +\infty }\|u(s,\cdot)\|_2^2=\begin{cases}
    +\infty\quad &\beta<2+\frac{8}{N},\\
    b^{\frac{N}{4-N}}\left(\|\nabla V\|_2^2\right)^{\frac{N}{4-N}} \|V\|_2^2& \beta=2+\frac{8}{N},\\
    0& \beta>2+\frac{8}{N}.
    \end{cases}
$$
\end{itemize}
\el

\bpr
(i) For any $s\in \R$, since $(\lambda(s), v(s,x))$ is a positive solution to \eqref{eq:20211029-wbe1}, by the definition of $t_s$ and $u(s,x)$, {\bf Theorem A} implies that $(\lambda(s), u(s,x))$ is a positive solution to \eqref{1.1}. Next, we shall prove that $\mathcal{C}$ is connected in $\R\times H^1(\R^N)$. Indeed, for any $s^*\in \R$ and $s_n\rightarrow s^*$, it follows that $\lambda(s_n)\rightarrow \lambda(s^*)$ since $\tilde{\mathcal{S}}$ is connected. Furthermore, $v(s_n, x)\rightarrow v(s^*,x)$ in $H^1(\R^N)$. So we have $\|\nabla v(s_n, \cdot)\|_2^2\rightarrow \|\nabla v(s^*, \cdot)\|_2^2$ as $n\rightarrow +\infty$.
And it follows by Lemma \ref{lemma:20211109-l1} that
$t_{s_n}\rightarrow t_{s^*}$ as $n\rightarrow +\infty$. Combining with $v(s_n,x)\rightarrow v(s^*, x)$ in $H^1(\R^N)$, we obtain that
$u(s_n,x)\rightarrow u(s^*, x)$ in $H^1(\R^N)$. Hence, $\mathcal{C}$ is connected in $\R\times H^1(\R^N)$.

(ii)
By \cite[Theorem 1.3-(i)]{JeanZhangZhong2021}, we have that $\|\nabla v(s,\cdot)\|_2\rightarrow 0$ as $s\rightarrow -\infty$. So by Lemma \ref{lemma:20211109-l1}, we have that
$t_s\rightarrow \frac{1}{\sqrt{a}}$ as $s\rightarrow -\infty$.
A direct computation shows that
$\|\nabla u(s,\cdot)\|_2^2=t_{s}^{2-N}\|\nabla v(s,\cdot)\|_2^2$ and
$\|u(s,\cdot)\|_2^2=t_{s}^{-N}\|v(s,\cdot)\|_2^2$.
Then \eqref{eq:20211109-be1} and \eqref{eq:20211109-be2} can also follow by \cite[Theorem 1.3-(i)]{JeanZhangZhong2021}

(iii) By \cite[Theorem 1.3-(ii)]{JeanZhangZhong2021}, we have that $\|\nabla v(s,\cdot)\|_2\rightarrow +\infty$ as $s\rightarrow +\infty$. So Lemma \ref{lemma:20211109-l1} implies that
$$\frac{t_s}{\left(\frac{1}{b\|\nabla v(s,\cdot)\|_2^2}\right)^{\frac{1}{4-N}}}\rightarrow 1\;\hbox{as}~s\rightarrow +\infty.$$
Recalling from \cite[Theorem 4.10]{JeanZhangZhong2021} that
$$\phi(s,x):=\lambda(s)^{\frac{1}{2-\beta}}v\left(s,\left(\frac{x}{\sqrt{\lambda(s)}}\right)\right)\rightarrow V~\hbox{in}~H^1(\R^N)~\hbox{as}~s\rightarrow +\infty$$
and a direct calculation shows that
$$
 \begin{cases}
 \|\nabla v(s,\cdot)\|_2^2=\lambda(s)^{1+\frac{2}{\beta-2}-\frac{N}{2}} \|\nabla \phi(s,\cdot)\|_2^2,\\
 \|v(s,\cdot)\|_2^2=\lambda(s)^{\frac{2}{\beta-2}-\frac{N}{2}}\|\phi(s,\cdot)\|_2^2.
 \end{cases}$$
 Then, we can arrive at
 \begin{align*}
 \lim_{s\rightarrow +\infty} \|\nabla u(s,\cdot)\|_2^2=&\lim_{s\rightarrow +\infty} \left(b \|\nabla v(s,\cdot)\|_2^2\right)^{\frac{N-2}{4-N}}\|\nabla v(s,\cdot)\|_2^2\\
 =&\lim_{s\rightarrow +\infty} b^{\frac{N-2}{4-N}}\left(\|\nabla v(s,\cdot)\|_2^2\right)^{\frac{2}{4-N}}\\
 =&b^{\frac{N-2}{4-N}}\left(\|\nabla V\|_2^2\right)^{\frac{2}{4-N}} \lim_{s\rightarrow +\infty} \lambda(s)^{[1+\frac{2}{\beta-2}-\frac{N}{2}]\cdot \frac{2}{4-N}}\\
 =&+\infty
 \end{align*}
and
\begin{align*}
\lim_{s\rightarrow +\infty} \|u(s,\cdot)\|_2^2=&\lim_{s\rightarrow +\infty} t_{s}^{-N}\| v(s,\cdot)\|_2^2\\
=&\lim_{s\rightarrow +\infty}\left(\frac{1}{b\|\nabla v(s,\cdot)\|_2^2}\right)^{\frac{-N}{4-N}}\|v(s,\cdot)\|_2^2\\
=&b^{\frac{N}{4-N}}\lim_{s\rightarrow +\infty} \left(\|\nabla v(s,\cdot)\|_2^2\right)^{\frac{N}{4-N}} \|v(s,\cdot)\|_2^2\\
=&b^{\frac{N}{4-N}}\left(\|\nabla V\|_2^2\right)^{\frac{N}{4-N}} \|V\|_2^2 \lim_{s\rightarrow +\infty} \lambda(s)^{[1+\frac{2}{\beta-2}-\frac{N}{2}]\cdot \frac{N}{4-N} +\frac{2}{\beta-2}-\frac{N}{2}}\\
=:&b^{\frac{N}{4-N}}\left(\|\nabla V\|_2^2\right)^{\frac{N}{4-N}} \|V\|_2^2 \lim_{s\rightarrow +\infty} \lambda(s)^\sigma,
\end{align*}
where
$$\sigma:=[1+\frac{2}{\beta-2}-\frac{N}{2}]\cdot \frac{N}{4-N} +\frac{2}{\beta-2}-\frac{N}{2}=\frac{8-N(\beta-2)}{(\beta-2)(4-N)}.$$
We note that
$$\sigma\begin{cases}
>0~\quad &\hbox{if}~\beta<2+\frac{8}{N},\\
=0&\hbox{if}~\beta=2+\frac{8}{N},\\
<0&\hbox{if}~\beta>2+\frac{8}{N}.
\end{cases}$$
Hence,
$$\lim_{s\rightarrow +\infty} \|u(s,\cdot)\|_2^2=\begin{cases}
+\infty~\quad &\hbox{if}~\beta<2+\frac{8}{N},\\
b^{\frac{N}{4-N}}\left(\|\nabla V\|_2^2\right)^{\frac{N}{4-N}} \|V\|_2^2 &\hbox{if}~\beta=2+\frac{8}{N},\\
0&\hbox{if}~\beta>2+\frac{8}{N}.
\end{cases}$$
\epr

\begin{figure}[htbp]
\centering
\begin{minipage}{0.8\linewidth}
\centering
\includegraphics[width=1\linewidth]{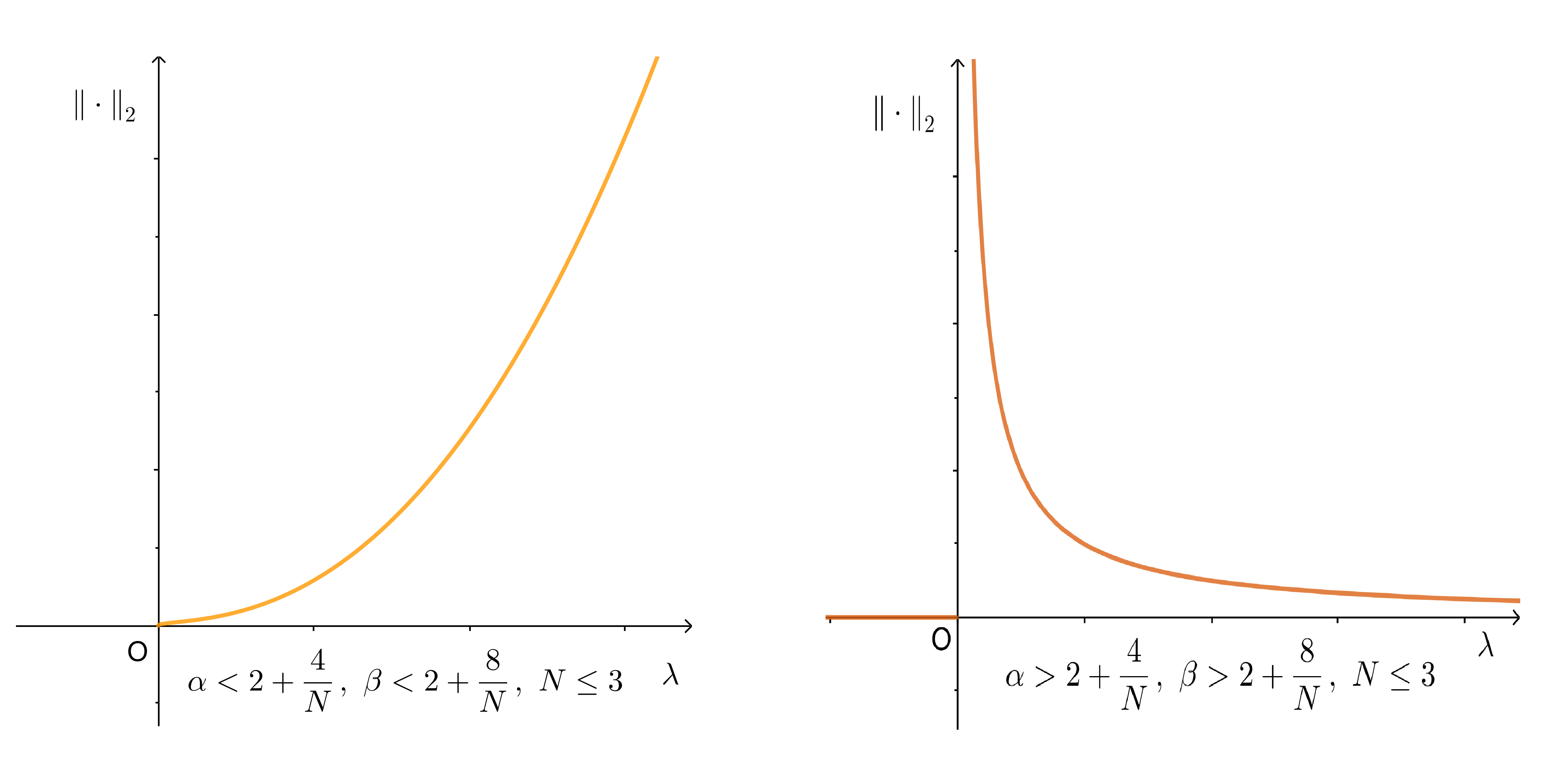}
\caption{The case of $N\leq 3$}
\label{Fig.1}
\end{minipage}
\begin{minipage}{0.8\linewidth}
\centering
\includegraphics[width=1\linewidth]{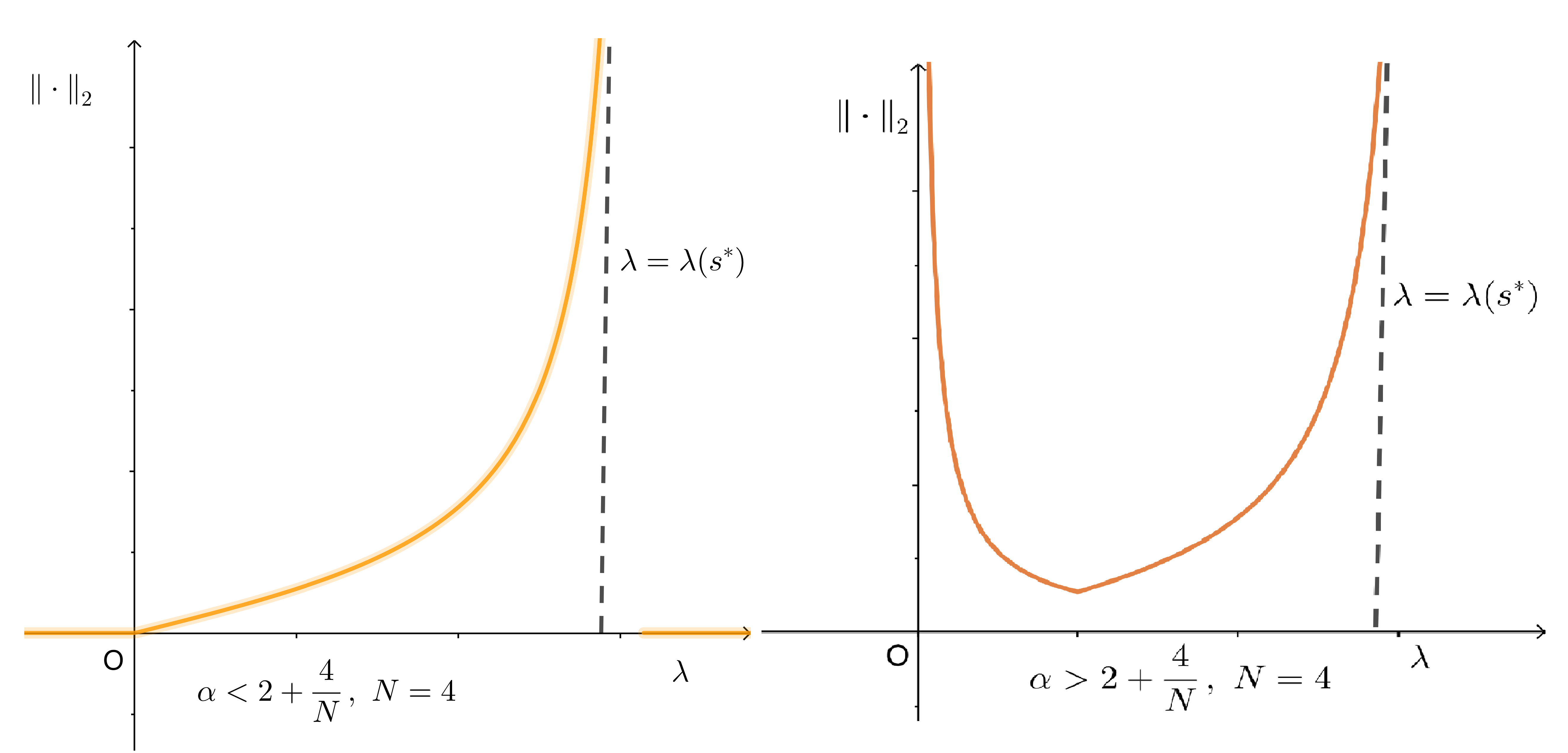}
\caption{The case of $N=4$}
\label{Fig.2}
\end{minipage}
\begin{minipage}{0.8\linewidth}
\centering
\includegraphics[width=1\linewidth]{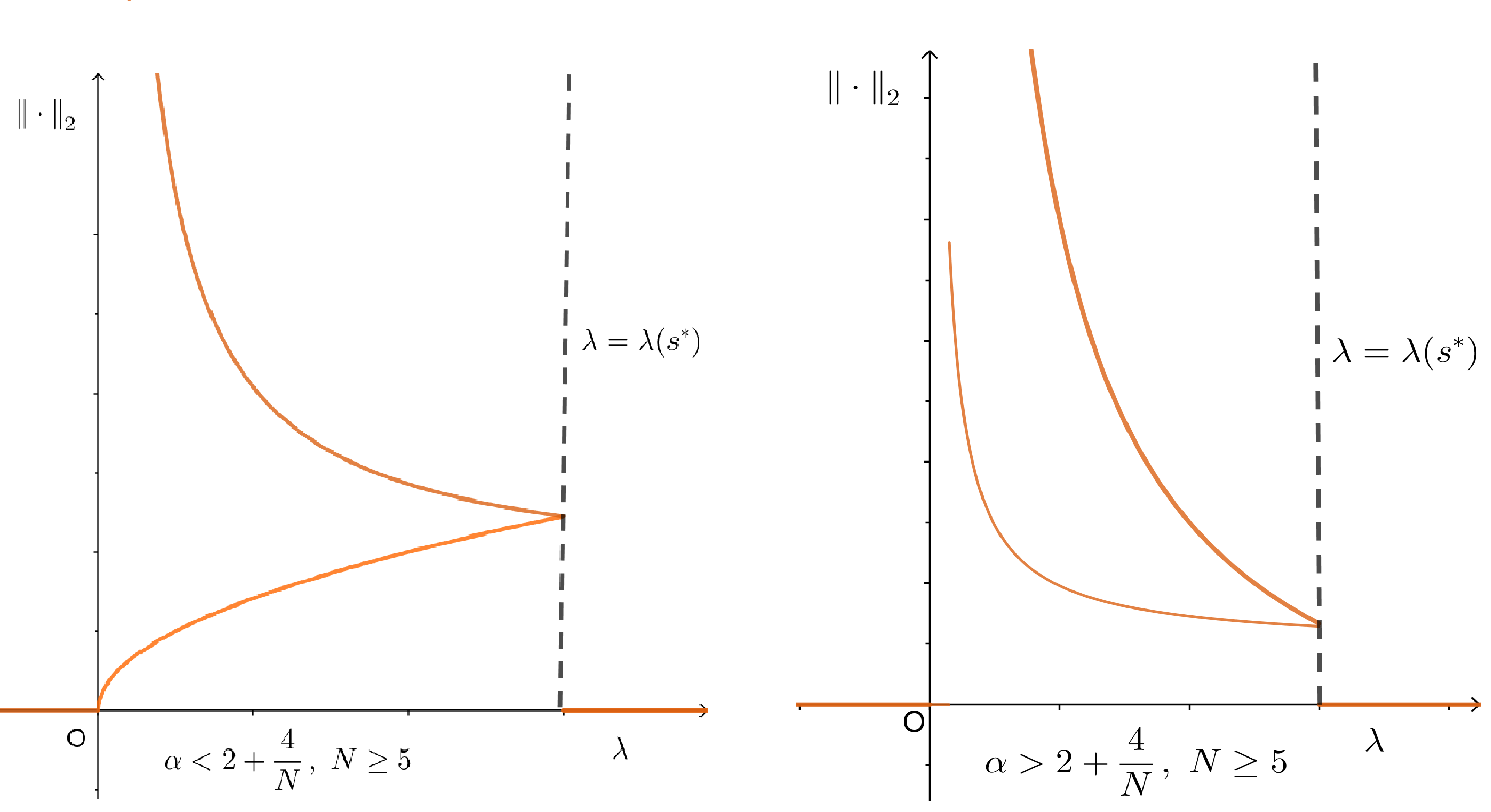}
\caption{The case of $N\geq 5$}
\label{Fig.3}
\end{minipage}
\end{figure}

Basing on Lemma \ref{lemma:20211109-l3}, we are clear with the $\|u\|_2^2$ along the curve $\mathcal{C}$. A schematic diagram for the cases of $K_{11}$ and $K_{33}$ can be seen in FIGURE\ref{Fig.1}.

\bl\lab{lemma:20211109-l4}
Assume that $a>0,b>0$ and $N\leq 3$. For any $0<\Lambda_1\leq \Lambda_2<+\infty$, define the set
$$
\mathcal{W}_{\Lambda_1}^{\Lambda_2}:=\left\{u\in H^1_{rad}(\R^N): \hbox{$u $ is a non-negative solution to \eqref{1.1} with $\lambda\in [\Lambda_1,\Lambda_2]$}\right\}.
$$
Then $\mathcal{W}_{\Lambda_1}^{\Lambda_2}$ is compact in $H^1(\R^N)$.
\el
\bpr
We also put
\beq\lab{eq:20211202-e3}
\mathcal{U}_{\Lambda_1}^{\Lambda_2}:=\left\{v\in H^1_{rad}(\R^N): \hbox{$v $ is a non-negative solution to \eqref{eq:20211029-wbe1} with $\lambda\in [\Lambda_1,\Lambda_2]$}\right\}.
\eeq
By {\bf Theorem A}, we see that $u\in \mathcal{W}_{\Lambda_1}^{\Lambda_2}$ if and only if
$u(x)=v(tx)$ with $v\in \mathcal{U}_{\Lambda_1}^{\Lambda_2}$ and $t=t_v$.
So for any sequence $\{u_n\}\subset \mathcal{W}_{\Lambda_1}^{\Lambda_2}$, there exists a corresponding sequence $\{v_n\}\subset \mathcal{U}_{\Lambda_1}^{\Lambda_2}$ such that
$u_n=v_n(t_n x)~\hbox{and}~t_n=t_{v_n}$.
So by $\mathcal{U}_{\Lambda_1}^{\Lambda_2}$ is compact (see \cite[Corollary 3.2]{JeanZhangZhong2021}), up to a subsequence, we assume that $v_n\rightarrow v\in \mathcal{U}_{\Lambda_1}^{\Lambda_2}$ in $H^1(\R^N)$. Put $u(x)=v(t_0 x)$ with $t_0:=t_v>0$, then $u\in \mathcal{W}_{\Lambda_1}^{\Lambda_2}$.
By Lemma \ref{lemma:20211109-l1}, $t_n\rightarrow t_0$ as $n\rightarrow+\infty$ and thus
$u_n=v_n(t_nx)\rightarrow v(t_0x)=u(x)$ in $H^1(\R^N)$.
Hence, $\mathcal{W}_{\Lambda_1}^{\Lambda_2}$ is compact in $H^1(\R^N)$.
\epr

\noindent
{\it Proof of Theorem \ref{th:20211109-th1}.} Basing on the Lemmas \ref{lemma:20211109-l1}-\ref{lemma:20211109-l4}, apply a similar argument as the proof of \cite[Theorem 1.1]{JeanZhangZhong2021}, one can obtain the conclusions of Theorem \ref{th:20211109-th1}. Here we only prove (1) and (4) as examples. Let us introduce the function
$$
  \rho: \mathcal{C}\rightarrow \mathbb{R}^+,\quad (\lambda(s),u(s,\cdot))\mapsto \|u(s,\cdot)\|_2^2.
$$

\noindent
(1)
By Lemma \ref{lemma:20211109-l3}-(i), there exists a connected positive solution branch of \eqref{1.1}, say $\mathcal{C}:=\left\{(\lambda(s), u(s,x)): s\in \R\right\}$. For $(\alpha,\beta)\in K_{11}$, by Lemma \ref{lemma:20211109-l3}-(ii), there exist $\{s_n\}\subset \R$ with $s_n\rightarrow -\infty$ and $\|u(s_n,\cdot)\|_2^2\rightarrow 0$.
Similarly, by Lemma \ref{lemma:20211109-l3}-(iii), there exist $\{s'_n\}\subset \R$ with $s'_n\rightarrow +\infty$ and $\|u(s'_n,\cdot)\|_2^2\rightarrow +\infty$.
Since $\mathcal{C}$ is connected, it follows that $\rho$ is onto.\\
(4) Let $\mathcal{S}:=\{(\lambda, u)\in (0, + \infty) \times H_{rad}^{1}(\RN)  : (\lambda, u)\;\hbox{solves \eqref{1.1}}, \, u>0 \}$
and $\tilde{\mathcal{S}}$ be the connected component containing $\mathcal{C}$.
So we have
$$\rho(\mathcal{S})\supset \rho(\widetilde{\mathcal{S}})\supset \rho(\mathcal{C})\supset\left(0,  b^{\frac{N}{4-N}}\left(\|\nabla V\|_2^2\right)^{\frac{N}{4-N}} \|V\|_2^2\right).$$
Hence, for any $a\in (0,b^{\frac{N}{4-N}}\left(\|\nabla V\|_2^2\right)^{\frac{N}{4-N}} \|V\|_2^2)$, \eqref{1.1} possesses at least one normalized solution $(\lambda, u_\lambda)$ with $\lambda>0$ and $0<u_\lambda\in H_{rad}^{1}(\RN)$.
On the other hand, by the local uniqueness result (see \cite[Theorem 5.1]{JeanZhangZhong2021}), we can find some $0<\Lambda_1<\Lambda_2<+\infty$ such that
$$P_2(\mathcal{S}\backslash \mathcal{C})\subset \mathcal{W}_{\Lambda_1}^{\Lambda_2},$$
where $P_2: (0, +\infty) \times  H_{rad}^{1}(\RN)  \rightarrow H_{rad}^{1}(\RN)$ be the projection onto the $H_{rad}^{1}(\RN)$-component.
By Lemma \ref{lemma:20211109-l4}, $\mathcal{W}_{\Lambda_1}^{\Lambda_2}$ is compact in $H^1(\R^N)$. It is well defined that
\begin{align*}
+\infty>c^*:=&\sup\left\{\rho(\lambda, u): (\lambda, u)\in \mathcal{S}\right\}\\
=&\max\Big\{\sup\big\{\|u\|_2^2:(\lambda,u)\in \mathcal{C}\big\}, \max\big\{\|u\|_2^2:u\in \mathcal{W}_{\Lambda_1}^{\Lambda_2}\big\}\Big\}\\
\geq& b^{\frac{N}{4-N}}\left(\|\nabla V\|_2^2\right)^{\frac{N}{4-N}} \|V\|_2^2,
\end{align*}
then for any $c>c^*$, \eqref{1.1} has no positive normalized solution.

More details we refer to \cite[Section 8]{JeanZhangZhong2021}.
\qed

\section{The case of $N=4$}\lab{sec:20230217-s2}
\renewcommand{\theequation}{3.\arabic{equation}}
\br\lab{remark:20211202-r1}
For $N=4$, it is easy to see that $at^2+b\|\nabla v\|_2^2-1=0$ has a positive root, if and only if $\|\nabla v\|_2^2<\frac{1}{b}$. In such a situation, the root $t_v$ is unique and $t_v=\left(\frac{1-b\|\nabla v\|_2^2}{a}\right)^{\frac{1}{2}}$.
\er

\noindent{\it Proof of Theorem \ref{thm:20211115-th1}.} Let $v(s,x)\in H^1(\R^N)$ be given in \eqref{eq:20230217-e1}.
We note that for $s$ sufficiently negative, $\|\nabla v(s,\cdot)\|_2^2$ is small, and then by Remark \ref{remark:20211202-r1}, there exists a unique $t_s>0$ such that
$$at_s^2+b\|\nabla v(s,\cdot)\|_2^2-1=0.$$
Recalling  \cite[Theorem 1.3]{JeanZhangZhong2021} that $\|\nabla v(s,\cdot)\|_2^2\rightarrow +\infty$ as $s\rightarrow +\infty$, there exists some $s^*\in \R$ such that
$$\|\nabla v(s,\cdot)\|_2^2<\frac{1}{b}~\hbox{for}~s\in (-\infty, s^*)~\hbox{and}~\|\nabla v(s^*,\cdot)\|_2^2=\frac{1}{b}.$$
So $t_s$ is well defined for $s\in (-\infty, s^*)$, furthermore, by Remark \ref{remark:20211202-r1} again,
$$\lim_{s\rightarrow -\infty}t_s=a^{-\frac{1}{2}}~\hbox{and}~\lim_{s\nearrow s^*}t_s=0.$$
So we also have that \eqref{eq:20211109-be2} holds. On the other hand,
\begin{align*}
\lim_{s\nearrow s^*}\|u(s,\cdot)\|_2^2=\lim_{s\nearrow s^*}t_{s}^{-N}\|v(s,\cdot)\|_2^2
=\|v(s^*,\cdot)\|_2^2 \lim_{s\nearrow s^*}t_{s}^{-N}
=+\infty.
\end{align*}
Hence, by the continuity, we conclude the final results. A schematic diagram for the cases of $\alpha<2+\frac{4}{N}$ and $\alpha>2+\frac{4}{N}$ with $N=4$ can be seen in FIGURE\ref{Fig.2}.\qed

\section{The case of $N\geq 5$}\lab{sec:20230217-s3}
\renewcommand{\theequation}{4.\arabic{equation}}
\bl\lab{lemma:20211115-l1}
Let $N\geq 5$.
For any $v\in H^1(\R^N)$, we put
$\displaystyle f(t):=at^2+b\|\nabla v\|_2^2t^{4-N}-1$ and
$t^*:=\left(\frac{b(N-4)\|\nabla v\|_2^2}{2a}\right)^{\frac{1}{N-2}}$, then we have that
$f(t)$ decreases in $(0,t^*)$ and $f(t)$ increases in $(t^*,+\infty)$.
\el
\bpr
By a direct computation, we consider
$$f'(t)=2at-b(N-4)\|\nabla v\|_2^2 t^{3-N}=0.$$
We see that $t=t^*$ is the unique root of $f'(t)=0$. In particular, for $t\in (0,t^*)$, we have that $f'(t)<0$ while $f'(t)>0$ for $t\in (t^*,+\infty)$.
\epr

\bl\lab{lemma:20211115-l2}
Let $v(s,x)\in H^1(\R^N)$ be given in \eqref{eq:20230217-e1}.
There exists some $s^*\in \R$ such that for any $s<s^*$,
there exist exactly two number $0<t_{1,s}<t^*_s<t_{2,s}<+\infty$ such that
$$at_{i,s}^{2}+b\|\nabla v(s,\cdot)\|_2^2t_{i,s}^{4-N}-1=0, i=1,2, s<s^*,$$
where
\beq\lab{eq:20211203-xe1}
t^*_s:=\left(\frac{b(N-4)\|\nabla v(s,\cdot)\|_2^2}{2a}\right)^{\frac{1}{N-2}}.
\eeq
Furthermore,
\beq\lab{eq:20211115-ze2}
\lim_{s\nearrow s^*}t_{i,s}=t_s^*, i=1,2,
\eeq
and
$$
\lim_{s\rightarrow 0^+}t_{i,s}=\begin{cases}
0\quad &\hbox{for}~i=1,\\
a^{-\frac{1}{2}}&\hbox{for}~i=2.
\end{cases}
$$
\el
\bpr
A direct corollary of Lemma \ref{lemma:20211115-l1} is that
$$
\min_{t>0}f(t)=f(t^*)=2^{-\frac{2}{N-2}}(N-2)(N-4)^{-\frac{N-4}{N-2}} a^{\frac{N-4}{N-2}} b^{\frac{2}{N-2}} \|\nabla v\|_{2}^{\frac{4}{N-2}}-1.
$$

Recalling that $\|\nabla v(s,\cdot)\|_2^2\rightarrow 0$ as  $s\rightarrow -\infty$, while $\|\nabla v(s,\cdot)\|_2^2\rightarrow +\infty$ as  $s\rightarrow +\infty$.
One can see that there exists some $s^*\in \R$ such that
\beq\lab{eq:20211202-e1}
2^{-\frac{2}{N-2}}(N-2)(N-4)^{-\frac{N-4}{N-2}} a^{\frac{N-4}{N-2}} b^{\frac{2}{N-2}} \|\nabla v(s^*,\cdot)\|_{2}^{\frac{4}{N-2}}-1=0
\eeq
and
$$
2^{-\frac{2}{N-2}}(N-2)(N-4)^{-\frac{N-4}{N-2}} a^{\frac{N-4}{N-2}} b^{\frac{2}{N-2}} \|\nabla v(s,\cdot)\|_{2}^{\frac{4}{N-2}}-1<0, \forall s<s^*.
$$
Now, for any $s<s^*$, let $t_s^*$ be given by \eqref{eq:20211203-xe1}.
Then  by Lemma \ref{lemma:20211115-l1}, we have that
$$\min_{t>0}[at^2+b\|\nabla v(s,\cdot)\|_2^2t^{4-N}-1]=a(t^*_s)^2+b\|\nabla v(s,\cdot)\|_2^2(t^*_s)^{4-N}-1<0.$$
And there exists a unique $t_{1,s}\in (0,t^*_s)$ and a unique $t_{2,s}\in (t^*_s,+\infty)$ such that
$$at_{i,s}^{2}+b\|\nabla v(s,\cdot)\|_2^2t_{i,s}^{4-N}-1=0, i=1,2.$$
We remark that the uniqueness implies that $t_{i,s}, i=1,2,$ are continuous in $(-\infty, s^*)$ (by applying a similar argument as Lemma \ref{lemma:20211109-l1}).
For $s\nearrow s^*$, recalling \eqref{eq:20211202-e1},  $a(t^*_s)^2+b\|\nabla v(s,\cdot)\|_2^2(t^*_s)^{4-N}-1=0$, one can see that \eqref{eq:20211115-ze2} holds. Furthermore, by $\|\nabla v(s,\cdot)\|_2^2\rightarrow 0$ as $s\rightarrow -\infty$, we see that $t^*_s\rightarrow 0$ as $s\rightarrow -\infty$. So by
$t_{1,s}\in (0, t_s^*)$ for $s\in (-\infty,s^*)$, we obtain that
$\displaystyle\lim_{s\rightarrow -\infty}t_{1,s}=0$.
On the other hand, by $t_{2,s}>t_s^*$ for $s\in (-\infty,s^*)$ and $N\geq 5$, we have that
\begin{align*}b\|\nabla v(s,\cdot)\|_2^2 t_{2,s}^{4-N}&\leq b\|\nabla v(s,\cdot)\|_2^2 (t_s^*)^{4-N}
\\&=(\frac{b(N-4)}{2a})^{\frac{4-N}{N-2}} \|\nabla v(s,\cdot)\|_{2}^{\frac{4}{N-2}}\rightarrow 0~\hbox{as}~s\rightarrow -\infty.\end{align*}
Then combining with $f(t_{2,s})=0$, we conclude that
$\displaystyle\lim_{s\rightarrow -\infty}t_{2,s}=a^{-\frac{1}{2}}$.
\epr

\br\lab{remark:20211115-r1}
For any $s<s^*$, we put
$$\mathcal{C}_1:=\{\big(\lambda(s),v(s, t_{1,s}x)\big): s<s^*\}~\hbox{and}~\mathcal{C}_2:=\{\big(\lambda(s),v(s, t_{2,s}x)\big): s<s^*\}.$$
Then we see that both $\mathcal{C}_1$ and $\mathcal{C}_2$ are curves in $\R\times H^1(\R^N)$.
Noting that for $s=s^*$, $u(x):=v(s^*, t_s^* x)$ is a solution to \eqref{1.1} with $\lambda=\lambda(s^*)$. We see that
$$\mathcal{C}:=\mathcal{C}_1\cup \{(\lambda(s^*),v(s^*, t_s^* x))\} \cup \mathcal{C}_2$$
is connected in $\R\times H^1(\R^N)$.
\er

\noindent
{\it Proof of Theorem \ref{thm:20211115-th2}.}
By $t_{2,s}\rightarrow a^{-\frac{1}{2}}$ as $s\rightarrow -\infty$, similar as before, we obtain that
$$
\lim_{s \rightarrow -\infty }\|u_2(s,\cdot)\|_2^2=\begin{cases}
    0\quad &\alpha<2+\frac{4}{N},\\
    a^{\frac{N}{2}}\|U\|_2^2& \alpha=2+\frac{4}{N},\\
    +\infty& \alpha>2+\frac{4}{N}.
    \end{cases}
$$
And by $t_{i,s}\rightarrow t_s^*$ as $s\rightarrow s^*, i=1,2$, we have
$$
\lim_{s \rightarrow s^* }\|u_i(s,\cdot)\|_2^2=\left(\frac{b(N-4)}{2a}\right)^{-\frac{N}{N-2}} \|\nabla v(s^*, \cdot)\|_{2}^{-\frac{2N}{N-2}} \|v(s^*, \cdot)\|_2^2, i=1,2.
$$
On the other hand, recalling from \cite[Theorem 4.6]{JeanZhangZhong2021} that
$$\psi(s,x):=\lambda(s)^{\frac{1}{2-\alpha}}v\left(s,\left(\frac{x}{\sqrt{\lambda(s)}}\right)\right)\rightarrow U~\hbox{in}~H^1(\R^N)~\hbox{as}~s\rightarrow -\infty,$$
then
\begin{align*}
&\hspace*{-8pt}t_{1,s}^{-N}\|v(s,\cdot)\|_2^2\\>&(t_s^*)^{-N}\|v(s,\cdot)\|_2^2\\
=&\left(\frac{b(N-4)}{2a}\right)^{-\frac{N}{N-2}}  \|\nabla v(s,\cdot)\|_{2}^{-\frac{2N}{N-2}} \cdot \|v(s,\cdot)\|_2^2\\
=&\left(\frac{b(N-4)}{2a}\right)^{-\frac{N}{N-2}}  \lambda(s)^{-[1+\frac{2}{\alpha-2}-\frac{N}{2}]\frac{N}{N-2}} \|\nabla \psi(s,\cdot)\|_{2}^{-\frac{2N}{N-2}} \lambda(s)^{\frac{2}{\alpha-2}-\frac{N}{2}} \|\psi(s,\cdot)\|_2^2\\
=&\left(\frac{b(N-4)}{2a}\right)^{-\frac{N}{N-2}} \|\nabla U\|_{2}^{-\frac{2N}{N-2}} \|U\|_2^2 \lambda(s)^{\eta} (1+o(1)),
\end{align*}
where
$$\eta:=\frac{2}{\alpha-2}-\frac{N}{2}-\left[1+\frac{2}{\alpha-2}-\frac{N}{2}\right]\frac{N}{N-2}=-\frac{4}{(N-2)(\alpha-2)}<0.$$
So by $\lambda(s)\rightarrow 0$ as $s\rightarrow -\infty$, we obtain that
$$\|u_{1}(s,\cdot)\|_2^2=t_{1,s}^{-N}\|v(s,\cdot)\|_2^2\rightarrow +\infty~\hbox{as}~s\rightarrow -\infty.$$
Hence, by Remark \ref{remark:20211115-r1}, we conclude the results by continuity.
In particular, we remark that
$$\frac{\|u_{1}(s,\cdot)\|_2^2}{\|u_{2}(s,\cdot)\|_2^2}\geq C \frac{\lambda(s)^{\eta}}{\lambda(s)^{\frac{2}{\alpha-2}-\frac{N}{2}}}
=C \lambda(s)^{\frac{N}{2(N-2)(\alpha-2)}[(N-2)\alpha-2N]}\rightarrow +\infty~\hbox{as}~s\rightarrow -\infty.$$

A schematic diagram for the cases of $\alpha<2+\frac{4}{N}$ and $\alpha>2+\frac{4}{N}$ with $N\geq 5$, we refer to FIGURE\ref{Fig.3}.\qed

\section{Proof of Theorem \ref{thm:20211115-th3}}\lab{sec:20230217-s4}

{\it Proof of Theorem \ref{thm:20211115-th3}.}
We argue by a way of negation. Suppose that there exists a sequence of positive solutions $(\lambda_n, u_n)$ such that $\|u_n\|_2^2:=c_n\rightarrow 0$ as $n\rightarrow +\infty$.
We put $t_n:=\frac{1}{\sqrt{a+\|\nabla u_n\|_2^2}}$ and $v_n(x):=u_n(\frac{x}{t_n})$, then we see that
$$-\Delta v_n+\lambda_n v_n=g(v_n)~\hbox{in}~\R^N.$$
We note that $t_n<\frac{1}{\sqrt{a}}, \forall n\in \N$, we obtain that
$$\|v_n\|_2^2=t_n^N\|u_n\|_2^2\rightarrow 0~\hbox{as}~n\rightarrow +\infty.$$
Since $\alpha\geq 2+\frac{4}{N}$, by \cite[Theorem 1.3]{JeanZhangZhong2021}, we see that $\{\lambda_n\}$ is bounded away from $0$, saying $\lambda_n\geq \Lambda_1>0$.
On the other hand, for any $\Lambda_2>\Lambda_1$, the set $\mathcal{U}_{\Lambda_1}^{\Lambda_2}$ defined by \eqref{eq:20211202-e3}, is compact in $H^1(\R^N)$ and bounded away from $0$, see also \cite[Lemma 3.3]{JeanZhangZhong2021}.
We conclude that $\lambda_n\rightarrow +\infty$.

Then by \cite[Theorem 1.3]{JeanZhangZhong2021} again, $\|\nabla v_n\|_2^2\rightarrow +\infty$ as $n\rightarrow +\infty$.
In particular, $t_n$ satisfies
$a t_n^2+b\|\nabla v_n\|_2^2 t_{n}^{4-N}=1$.
However, for $N=4$, we have that
$1=a t_n^2+b\|\nabla v_n\|_2^2 t_{n}^{4-N}>b\|\nabla v_n\|_2^2\rightarrow +\infty$,
which is a contradiction.
And for the case of $N\geq 5$,
\begin{align*}
1=&a t_n^2+b\|\nabla v_n\|_2^2 t_{n}^{4-N}\geq \min_{t>0} \left[at^2+b\|\nabla v_n\|_2^2 t^{4-N}\right]\\
=&2^{-\frac{2}{N-2}}(N-2)(N-4)^{-\frac{N-4}{N-2}} a^{\frac{N-4}{N-2}} b^{\frac{2}{N-2}} \|\nabla v_n\|_{2}^{\frac{4}{N-2}}
\rightarrow +\infty~\hbox{as}~n\rightarrow +\infty,
\end{align*}
also yields a contradiction.
\qed

\section*{Acknowledgements}
The authors would like to express their sincere gratitude to the anonymous referee for her/his valuable suggestions and comments.


\begin{thebibliography}{99}
\bibitem{Alves2005} (MR2123187) [10.1016/j.camwa.2005.01.008]
\newblock C.~O.~Alves, F.~J.~S.~A.~Corr\^{e}a and T.~F.~Ma,
\newblock {Positive solutions for a quasilinear elliptic equation of {K}irchhoff type},
\newblock \emph{Comput. Math. Appl.}, \textbf{49} (2005), 85-93.

\bibitem{Azzollini2012} (MR2951740)
\newblock A.~Azzollini,
\newblock {The elliptic {K}irchhoff equation in {$\Bbb R^N$} perturbed by a local nonlinearity},
\newblock \emph{Differential Integral Equations}, \textbf{25} (2012), 543-554.

\bibitem{Azzollini2015} (MR3359229) [10.1142/S0219199714500394]
\newblock A.~Azzollini,
\newblock { A note on the elliptic {K}irchhoff equation in {$\Bbb{R}^N$} perturbed by a local nonlinearity},
\newblock \emph{Commun. Contemp. Math.}, \textbf{17} (2015), 1450039, 5 pp.

\bibitem{BartschValeriola2013} (MR3009665) [10.1007/s00013-012-0468-x]
\newblock T.~Bartsch and S.~de~Valeriola,
\newblock {Normalized solutions of nonlinear {S}chr\"{o}dinger equations},
\newblock \emph{Arch. Math. (Basel)}, \textbf{100} (2013), 75-83.

\bibitem{BartschSoave2017} (MR3639521) [10.1016/j.jfa.2017.01.025]
\newblock T.~Bartsch and N.~Soave,
\newblock {A natural constraint approach to normalized solutions of nonlinear {S}chr\"{o}dinger equations and systems},
\newblock \emph{J. Funct. Anal.}, \textbf{272} (2017), 4998-5037.

\bibitem{Berestycki1983} (MR695535) [10.1007/BF00250555]
\newblock H.~Berestycki and P.-L.~Lions,
\newblock {Nonlinear scalar field equations. {I}. {E}xistence of a ground state},
\newblock \emph{Arch. Rational Mech. Anal.}, \textbf{82} (1983), 313-345.

\bibitem{Berestycki1983b} (MR734575)
\newblock H.~Berestycki, T.~Gallou\"{e}t and O.~Kavian,
\newblock {\'{E}quations de champs scalaires euclidiens non lin\'{e}aires dans le plan},
\newblock \emph{C. R. Acad. Sci. Paris S\'{e}r. I Math.}, \textbf{297} (1983), 307-310.

\bibitem{Berestycki1983a} (MR695536) [10.1007/BF00250556]
\newblock H.~Berestycki and P.-L.~Lions,
\newblock {Nonlinear scalar field equations. {II}. {E}xistence of infinitely many solutions},
\newblock \emph{Arch. Rational Mech. Anal.}, \textbf{82} (1983), 347-375.

\bibitem{BieganowskiMederski} (MR4232669) [10.1016/j.jfa.2021.108989]
\newblock B.~Bieganowski and J.~Mederski,
\newblock {Normalized ground states of the nonlinear {S}chr\"{o}dinger equation with at least mass critical growth},
\newblock \emph{J. Funct. Anal.}, \textbf{280} (2021), Paper No. 108989, 26 pp.

\bibitem{Chen2013} (MR3004514) [10.1016/j.nonrwa.2012.10.010]
\newblock S.~J.~Chen and L.~Li,
\newblock {Multiple solutions for the nonhomogeneous {K}irchhoff equation on {$\mathbb{R}^N$}},
\newblock \emph{Nonlinear Anal. Real World Appl.}, \textbf{14} (2013), 1477-1486.

\bibitem{Deng2015} (MR3406860) [10.1016/j.jfa.2015.09.012]
\newblock Y.~Deng, S.~Peng and W.~Shuai,
\newblock {Existence and asymptotic behavior of nodal solutions for the {K}irchhoff-type problems in {$\mathbb{R}^3$}},
\newblock \emph{J. Funct. Anal.}, \textbf{269} (2015), 3500-3527.

\bibitem{Figueiredo2014} (MR3218834) [10.1007/s00205-014-0747-8]
\newblock G.~M.~Figueiredo, N.~Ikoma and J.~R. Santos~J\'{u}nior,
\newblock {Existence and concentration result for the {K}irchhoff type equations with general nonlinearities},
\newblock \emph{Arch. Ration. Mech. Anal.}, \textbf{213} (2014), 931-979.

\bibitem{Guo2018} (MR3809133) [10.3934/cpaa.2018089]
\newblock H.~Guo, Y.~Zhang and H.~Zhou,
\newblock {Blow-up solutions for a {K}irchhoff type elliptic equation with trapping potential},
\newblock \emph{Commun. Pure Appl. Anal.}, \textbf{17} (2018), 1875-1897.

\bibitem{Guo2015} (MR3360660) [10.1016/j.jde.2015.04.005]
\newblock Z.~Guo,
\newblock {Ground states for {K}irchhoff equations without compact condition},
\newblock \emph{J. Differential Equations}, \textbf{259} (2015), 2884-2902.

\bibitem{HeLvZhangZhong2021} (MR4546673) [10.1016/j.jde.2023.01.039]
\newblock Q.~He, Z.~Lv, Y.~Zhang and X.~Zhong,
\newblock {Existence and blow up behavior of positive normalized solution to the Kirchhoff equation with general nonlinearities: Mass super-critical case}
\newblock \emph{J. Differential Equations}, \textbf{356} (2023), 375-406.

\bibitem{He2012} (MR2853562) [10.1016/j.jde.2011.08.035]
\newblock X.~He and W.~Zou,
\newblock {Existence and concentration behavior of positive solutions for a {K}irchhoff equation in {$\Bbb R^3$}},
\newblock \emph{J. Differential Equations}, \textbf{252} (2012), 1813-1834.

\bibitem{He2016} (MR3552562) [10.1016/j.jde.2016.08.034]
\newblock Y.~He,
\newblock {Concentrating bounded states for a class of singularly perturbed {K}irchhoff type equations with a general nonlinearity}
\newblock \emph{J. Differential Equations}, \textbf{261} (2016), 6178-6220.

\bibitem{He2014} (MR3194366) [10.1515/ans-2014-0214]
\newblock Y.~He, G.~Li and S.~Peng,
\newblock {Concentrating bound states for {K}irchhoff type problems in {$\Bbb R^3$} involving critical {S}obolev exponents},
\newblock \emph{Adv. Nonlinear Stud.}, \textbf{14} (2014), 483-510.

\bibitem{ikomatanaka2019} (MR4021261)
\newblock N.~Ikoma and K.~Tanaka,
\newblock {A note on deformation argument for {$L^2$} normalized solutions of nonlinear {S}chr\"{o}dinger equations and systems},
\newblock \emph{Adv. Differential Equations}, \textbf{24} (2019), 609-646.

\bibitem{Jeanjean1997} (MR1430506) [10.1016/S0362-546X(96)00021-1]
\newblock L.~Jeanjean,
\newblock {Existence of solutions with prescribed norm for semilinear elliptic equations},
\newblock \emph{Nonlinear Anal.}, \textbf{28} (1997), 1633-1659.

\bibitem{JeanjeanLe2022} (MR4476243) [10.1007/s00208-021-02228-0]
\newblock L.~Jeanjean and T.~T. Le,
\newblock {Multiple normalized solutions for a Sobolev critical {S}chr\"odinger equation},
\newblock \emph{Math. Ann.}, \textbf{384} (2022), 101-134.

\bibitem{JeanLu2019} (MR4030599) [10.1088/1361-6544/ab435e]
\newblock L.~Jeanjean and S.-S.~Lu,
\newblock {Nonradial normalized solutions for nonlinear scalar field equations},
\newblock \emph{Nonlinearity}, \textbf{32} (2019), 4942-4966.

\bibitem{Jeanjean-Lu-2022} (MR4493269) [10.1007/s00526-022-02320-6]
\newblock L.~Jeanjean and S.-S.~Lu,
\newblock {On global minimizers for a mass constrained problem},
\newblock \emph{Calc. Var. Partial Differential Equations}, \textbf{61} (2022), Paper No. 214, 18 pp.

\bibitem{JeanLu2020a} (MR4150876) [10.1007/s00526-020-01828-z]
\newblock L.~Jeanjean and S.-S.~Lu,
\newblock {A mass supercritical problem revisited},
\newblock \emph{Calc. Var. Partial Differential Equations}, \textbf{59} (2020), Paper No. 174, 43 pp.

\bibitem{Jeanjean2003} (MR2017241) [10.1515/ans-2003-0403]
\newblock L.~Jeanjean and K.~Tanaka,
\newblock {A note on a mountain pass characterization of least energy solutions},
\newblock \emph{Adv. Nonlinear Stud.}, \textbf{3} (2003), 445-455.

\bibitem{Jeanjean2003a} (MR1974637) [10.1090/S0002-9939-02-06821-1]
\newblock L.~Jeanjean and K.~Tanaka,
\newblock {A remark on least energy solutions in {$\mathbb{R}^N$}},
\newblock \emph{Proc. Amer. Math. Soc.}, \textbf{131} (2003), 2399-2408.

\bibitem{JeanZhangZhong2021}
\newblock L.~Jeanjean, J.~Zhang and X.~Zhong,
\newblock {A global branch approach to normalized solutions for the {S}chr\"odinger equation},
\newblock {arXiv e-prints}, (2021). {2112.05869}

\bibitem{Kirchhoff}
\newblock G.~Kirchhoff,
\newblock \emph{Mechanik, Teubner, Leipzig},
\newblock 1883.

\bibitem{Lei2015} (MR3250494) [10.1016/j.jmaa.2014.07.031]
\newblock C.-Y.~Lei, J.-F.~Liao and C.-L.~Tang,
\newblock {Multiple positive solutions for {K}irchhoff type of problems with singularity and critical exponents},
\newblock \emph{J. Math. Anal. Appl.}, \textbf{421} (2015), 521-538.

\bibitem{Li2021} (MR4443235) [10.54330/afm.120247]
\newblock G.~Li, X.~Luo and T.~Yang,
\newblock {Normalized solutions to a class of Kirchhoff equations with Sobolev critical exponent},
\newblock \emph{Ann. Fenn. Math.}, \textbf{47} (2022),  895-925.

\bibitem{Li2014} (MR3200382) [10.1016/j.jde.2014.04.011]
\newblock G.~Li and H.~Ye,
\newblock {Existence of positive ground state solutions for the nonlinear {K}irchhoff type equations in {$\Bbb R^3$}},
\newblock \emph{J. Differential Equations}, \textbf{257} (2014), 566-600.

\bibitem{Li2019} (MR3987389) [10.1016/j.na.2018.12.010]
\newblock Y.~Li, X.~Hao and J.~Shi,
\newblock {The existence of constrained minimizers for a class of nonlinear {K}irchhoff-{S}chr\"{o}dinger equations with doubly critical exponents in dimension four},
\newblock \emph{Nonlinear Anal.}, \textbf{186} (2019), 99-112.

\bibitem{Li2012} (MR2946973) [10.1016/j.jde.2012.05.017]
\newblock Y.~Li, F.~Li and J.~Shi,
\newblock {Existence of a positive solution to {K}irchhoff type problems without compactness conditions},
\newblock \emph{J. Differential Equations}, \textbf{253} (2012), 2285-2294.

\bibitem{Lions1978} (MR519648)
\newblock J.-L.~Lions,
\newblock {On some questions in boundary value problems of mathematical physics},
\newblock in \emph{Contemporary Developments in Continuum Mechanics and Partial Differential Equations ({P}roc. {I}nternat. {S}ympos., {I}nst. {M}at., {U}niv. {F}ed. {R}io de {J}aneiro, {R}io de {J}aneiro, 1977)}, volume~30 of \emph{North-Holland Math. Stud.}, pages 284-346. North-Holland, Amsterdam-New York, 1978.

\bibitem{Lions1984a} (MR778970) [10.1016/s0294-1449(16)30428-0]
\newblock P.-L.~Lions,
\newblock {The concentration-compactness principle in the calculus of variations. {T}he locally compact case. {I}},
\newblock \emph{Ann. Inst. H. Poincar\'{e} Anal. Non Lin\'{e}aire}, \textbf{1} (1984), 109-145.

\bibitem{Lions1984b} (MR778974) [10.1016/s0294-1449(16)30422-x]
\newblock P.-L.~Lions,
\newblock {The concentration-compactness principle in the calculus of variations. {T}he locally compact case. {II}},
\newblock \emph{Ann. Inst. H. Poincar\'{e} Anal. Non Lin\'{e}aire}, \textbf{1} (1984), 223-283.

\bibitem{Lu2017} (MR3567960) [10.1016/j.nonrwa.2016.09.003]
\newblock S.-S.~Lu,
\newblock {An autonomous {K}irchhoff-type equation with general nonlinearity in {$\Bbb{R}^N$}},
\newblock \emph{Nonlinear Anal. Real World Appl.}, \textbf{34} (2017), 233-249.

\bibitem{Luo2017} (MR3543109) [10.1016/j.nonrwa.2016.06.001]
\newblock X.~Luo and Q.~Wang,
\newblock {Existence and asymptotic behavior of high energy normalized solutions for the {K}irchhoff type equations in {$\Bbb{R}^3$}},
\newblock \emph{Nonlinear Anal. Real World Appl.}, \textbf{33} (2017), 19-32.

\bibitem{Meng2022} (MR4322298) [10.1016/j.jmaa.2021.125727]
\newblock X.~Meng and X.~Zeng,
\newblock {Existence and asymptotic behavior of minimizers for the {K}irchhoff functional with periodic potentials},
\newblock \emph{J. Math. Anal. Appl.}, \textbf{507} (2022), Paper No. 125727, 17 pp.

\bibitem{Naimen2014} (MR3278854) [10.1007/s00030-014-0271-4]
\newblock D.~Naimen,
\newblock {Positive solutions of {K}irchhoff type elliptic equations involving a critical {S}obolev exponent},
\newblock \emph{NoDEA Nonlinear Differential Equations Appl.}, \textbf{21} (2014), 885-914.

\bibitem{Shibata2014} (MR3147450) [10.1007/s00229-013-0627-9]
\newblock M.~Shibata,
\newblock {Stable standing waves of nonlinear {S}chr\"{o}dinger equations with a general nonlinear term},
\newblock \emph{Manuscripta Math.}, \textbf{143} (2014), 221-237.

\bibitem{Shuai2015} (MR3345850) [10.1016/j.jde.2015.02.040]
\newblock W.~Shuai,
\newblock {Sign-changing solutions for a class of {K}irchhoff-type problem in bounded domains},
\newblock \emph{J. Differential Equations}, \textbf{259} (2015), 1256-1274.

\bibitem{Soave-1} (MR4107073) [10.1016/j.jde.2020.05.016]
\newblock N.~Soave,
\newblock {Normalized ground states for the {NLS} equation with combined nonlinearities},
\newblock \emph{J. Differential Equations}, \textbf{269} (2020), 6941-6987.

\bibitem{Soave2020} (MR4096725) [10.1016/j.jfa.2020.108610]
\newblock N.~Soave,
\newblock {Normalized ground states for the {NLS} equation with combined nonlinearities: The {S}obolev critical case},
\newblock \emph{J. Funct. Anal.}, \textbf{279} (2020), 108610, 43 pp.

\bibitem{Stefanov19} (MR3975859) [10.1007/s00220-019-03484-7]
\newblock A.~Stefanov,
\newblock {On the normalized ground states of second order {PDE}'s with mixed power non-linearities},
\newblock \emph{Comm. Math. Phys.}, \textbf{369} (2019), 929-971.

\bibitem{Stuart1982} (MR662670) [10.1112/plms/s3-45.1.169]
\newblock C.~A.~Stuart,
\newblock {Bifurcation for {D}irichlet problems without eigenvalues},
\newblock \emph{Proc. London Math. Soc. (3)}, \textbf{45} (1982), 169-192.

\bibitem{Stuart1981} (MR819032)
\newblock A.~Stuart,
\newblock {Bifurcation from the continuous spectrum in the {$L^2$}-theory of elliptic equations on {$\mathbb{R}^n$}},
\newblock in \emph{Recent Methods in Nonlinear Analysis and Applications ({N}aples, 1980)}, pages 231-300. Liguori, Naples, (1981).

\bibitem{Wang2012} (MR2946975) [10.1016/j.jde.2012.05.023]
\newblock J.~Wang, L.~Tian, J.~Xu and F.~Zhang,
\newblock {Multiplicity and concentration of positive solutions for a {K}irchhoff type problem with critical growth},
\newblock \emph{J. Differential Equations}, \textbf{253} (2012), 2314-2351.

\bibitem{Wang2020} (MR4088927) [10.1002/mma.6256]
\newblock Z.~Wang, X.~Zeng and Y.~Zhang,
\newblock {Multi-peak solutions of {K}irchhoff equations involving subcritical or critical {S}obolev exponents},
\newblock \emph{Math. Methods Appl. Sci.}, \textbf{43} (2020), 5151-5161.

\bibitem{Wei-Wu2022} (MR4433054) [10.1016/j.jfa.2022.109574]
\newblock J.~Wei and Y.~Wu,
\newblock {Normalized solutions for {S}chr\"odinger equations with critical Sobolev exponent and mixed nonlinearities},
\newblock \emph{J. Funct. Anal.}, \textbf{283} (2022), Paper No. 109574, 46 pp.

\bibitem{Xie2022} (MR4338596) [10.1007/s00033-021-01626-3]
\newblock Q.~Xie and B.-X.~Zhou,
\newblock {A study on the critical Kirchhoff problem in high-dimensional space},
\newblock \emph{Z. Angew. Math. Phys.}, \textbf{73} (2022), Paper No. 4, 29 pp.

\bibitem{Ye2015a} (MR3377698) [10.1007/s00033-014-0474-x]
\newblock H.~Ye,
\newblock {The existence of normalized solutions for {$L^2$}-critical constrained problems related to {K}irchhoff equations},
\newblock \emph{Z. Angew. Math. Phys.}, \textbf{66} (2015), 1483-1497.

\bibitem{Ye2015} (MR3382697) [10.1002/mma.3247]
\newblock H.~Ye,
\newblock {The sharp existence of constrained minimizers for a class of nonlinear {K}irchhoff equations},
\newblock \emph{Math. Methods Appl. Sci.}, \textbf{38} (2015), 2663-2679.

\bibitem{Ye2016} (MR3483878) [10.1007/s00033-016-0624-4]
\newblock H.~Ye,
\newblock {The mass concentration phenomenon for {$L^2$}-critical constrained problems related to {K}irchhoff equations},
\newblock \emph{Z. Angew. Math. Phys.}, \textbf{67} (2016), Art. 29, 16 pp.

\bibitem{Zeng2017} (MR3677841) [10.1016/j.aml.2017.05.012]
\newblock X.~Zeng and Y.~Zhang,
\newblock {Existence and uniqueness of normalized solutions for the {K}irchhoff equation},
\newblock \emph{Appl. Math. Lett.}, \textbf{74} (2017), 52-59.

\bibitem{Zhang2021} (MR4373844) [10.1063/5.0067520]
\newblock P.~Zhang and Z.~Han,
\newblock {Normalized ground states for Kirchhoff equations in ${\mathbb{R}}^{3}$ with a critical nonlinearity},
\newblock \emph{J. Math. Phys.}, \textbf{63} (2022),  Paper No. 021505, 15 pp.

\bibitem{Zhang2020} (MR4079057) [10.1016/j.na.2020.111856]
\newblock Y.~Zhang, X.~Tang and D.~Qin,
\newblock {Infinitely many solutions for {K}irchhoff problems with lack of compactness},
\newblock \emph{Nonlinear Anal.}, \textbf{197} (2020), 111856, 31 pp.

\bibitem{Zhu2021} (MR4311465) [10.1007/s00009-021-01835-0]
\newblock X.~Zhu, C.~Wang and Y.~Xue,
\newblock {Constraint minimizers of {K}irchhoff-{S}chr\"{o}dinger energy functionals with {$L^2$}-subcritical perturbation},
\newblock \emph{Mediterr. J. Math.}, \textbf{18} (2021), Paper No. 224, 20 pp.



\end{thebibliography}

\medskip
Received May 2023; revised July 2023; early access August 2023.
\medskip

\end{document}